\documentclass[10pt]{article}

\usepackage{graphicx}%
\usepackage{amsmath,amssymb,amsthm,upref,bm}%

\newtheorem{theorem}{Theorem}%
\newtheorem{remark}{Remark}

\begin{document}

\baselineskip=4.4mm

\makeatletter

\newcommand{\E}{\mathrm{e}\kern0.2pt} 
\newcommand{\D}{\mathrm{d}\kern0.2pt}
\newcommand{\RR}{\mathbb{R}}
\newcommand{\CC}{\mathbb{C}}%
\newcommand{\ii}{\kern0.05em\mathrm{i}\kern0.05em}

\renewcommand{\Re}{\mathrm{Re}} 
\renewcommand{\Im}{\mathrm{Im}}

\def\bottomfraction{0.9}

\title{\bf On characterization of balls via solutions \\ to the Helmholtz equation}

\author{Nikolay Kuznetsov}

\date{}

\maketitle

\vspace{-8mm}

\begin{center}
Laboratory for Mathematical Modelling of Wave Phenomena, \\ Institute for Problems
in Mechanical Engineering, Russian Academy of Sciences, \\ V.O., Bol'shoy pr. 61,
St Petersburg 199178, Russian Federation \\ E-mail address:
nikolay.g.kuznetsov@gmail.com
\end{center}

\begin{abstract}
\noindent A new analytical characterization of balls in the Euclidean space $\RR^m$
is obtained. Unlike previous results of this kind, using either harmonic functions
or solutions to the modified Helmholtz equation, the present one is based on
solutions to the Helmholtz equation. This is achieved at the expense of a
restriction imposed on the size of a domain---a feature absent in the inverse mean
value properties known before.
\end{abstract}

\setcounter{equation}{0}


\section{Introduction and notation}

In 1972 Kuran \cite{K} obtained the following characterization of balls via
harmonic functions:
\begin{quote}
{\it Let $D$ be a domain (= connected open set) of finite (Lebesgue) measure in the
Euclidean space ${\bf R}^m$ where $m \geq 2$. Suppose that there exists a point
$P_0$ in $D$ such that, for every function $h$ harmonic in $D$ and integrable over
$D$, the volume mean of $h$ over $D$ equals $h (P_0)$. Then $D$ is an open ball
(disk when $m=2$) centred at $P_0$.}
\end{quote}
A slight modification of Kuran's considerations shows that his theorem is valid even
if~$D$ is disconnected; see \cite{NV}, p.~377. This result was originally proved by
Epstein \cite{E} for a simply connected two-dimensional~$D$. Armitage and Goldstein
\cite{AG} reinforced this assertion, assuming that the mean value identity holds
only for positive harmonic functions which are $L^p$-integrable, $p \in  (0, n /
(n-2))$.

In the survey article \cite{NV}, one also finds a discussion of applications of
Kuran's theorem and of possible similar results involving some kinds of average
over $\partial D$, when $D$ is a bounded domain. A great deal of other interesting
material obtained in this area is reviewed in \cite{Ku4}; in particular, other
characterizations of balls via harmonic functions as well as of other domains
(strips, annuli etc.). Moreover, it contains a characterization of balls by
solutions to the modified Helmholtz equation
\begin{equation}
\nabla^2 u - \mu^2 u = 0 , \quad \mu \in \RR \setminus \{0\} . \label{MHh}
\end{equation}
Here and below, $\nabla = (\partial_1, \dots , \partial_m)$, $\partial_i = \partial
/ \partial x_i$, denotes the gradient operator. The latter result (discovered by
the author in 2021; see \cite{Ku}) is closely related to that obtained in this note.
Therefore, it is reasonable to formulate it here, but before that we introduce some
notation.

Let $x = (x_1, \dots, x_m)$ be a point in $\RR^m$, $m \geq 2$, by $B_r (x) = \{ y
\in \RR^m : |y-x| < r \}$ we denote the open ball of radius $r$ centred at $x$ (just
$B_r$, if centred at the origin). The ball is called admissible with respect to a
domain $D \subset \RR^m$ provided $\overline{B_r (x)} \subset D$. If $D$ has finite
Lebesgue measure and a function $f$ is integrable over $D$, then
\[ M (f, D) = \frac{1}{|D|} \int_{D} f (x) \, \D x
\]
is its volume mean value over $D$. Here and below $|D|$ is the domain's volume (area
if $D \subset \RR^2$); the volume of $B_r$ is $|B_r| = \omega_m r^m$, where
$\omega_m = 2 \, \pi^{m/2} / [m \Gamma (m/2)]$ is the volume of the unit ball; as
usual $\Gamma$ denotes the Gamma function. A dilated copy of a bounded domain $D$ is
$D_r = D \cup \left[ \cup_{x \in \partial D} B_r (x) \right]$. Thus, the distance
from $\partial D_r$ to $D$ is equal to $r$. Now we are in a position to formulate
the following assertion proved in~\cite{Ku}.

\begin{theorem}
Let $D \subset \RR^m$, $m \geq 2$, be a bounded domain, whose complement is
connected, and let $r > 0$ be such that $|B_r| = |D|$. If for a point $x_0 \in D$
and some $\mu > 0$ the identity
\begin{equation}
u (x_0) \, \Gamma \left(\frac{m}{2} + 1\right) \frac{I_{m/2} (\mu r)}{(\mu r /
2)^{m/2}} = M (u, D) \label{CR}
\end{equation}
holds for every positive function $u$ satisfying equation \eqref{MHh} in $D_r$, then
$D = B_r (x_0)$ $(I_\nu$ stands for the modified Bessel function of order $\nu)$.
\end{theorem}

In this note, we prove another analytic characterization of balls, which is also
based on the $m$-dimensional volume mean value identity, but instead of solutions
equation \eqref{MHh} it involves real-valued solutions of the Helmholtz equation
\begin{equation}
\nabla^2 u + \lambda^2 u = 0 \, , \label{Hh}
\end{equation}
here $\lambda$ is an arbitrary nonzero real number, say positive.

\section{Main result and discussion}

Prior to formulating the main result, we recall the $m$-dimensional mean value
formula:
\begin{equation}
a_m (\lambda r) \, u (x) = M (u, B_r (x)) \, , \ \ \mbox{where} \ \ a_m (t) = \Gamma
\left( \frac{m}{2} + 1 \right) \frac{J_{m/2} (t)}{(t / 2)^{m/2}} \, .
\label{MM}
\end{equation}
It holds for every admissible ball $B_r (x)$ provided $u \in C^2 (D) \cap L^1 (D)$
is a solution of \eqref{Hh} in $D$. As usual, $J_\nu$ denotes the Bessel function of
order $\nu$; its $n$th positive zero is denoted by $j_{\nu,n}$ (this standard
notation is used below). It is quite strange that equality \eqref{MM} was obtained
by the author only recently (see \cite{Ku1}, p.~84); before that only the spherical
mean value formula was known (see, for example, \cite{CH}, pp.~317--320). A detailed
discussion of identity \eqref{MM} can be found in~\cite{Ku2}, where the analogous
formula for solutions of \eqref{MHh} is also derived. Notice that $a_m (\lambda r)$
oscillates about the zero, whereas the similar expression in identity \eqref{CR} is
positive and increases monotonically with~$r$. The role played by monotonicity in
the proof of Theorem~1 suggests the following modification, which is our main
result.

\begin{theorem}
Let $D \subset \RR^m$, $m \geq 2$, be a bounded domain, whose complement is
connected, and let $r > 0$ be such that $|B_r| = |D|$. Suppose that for some
$\lambda > 0$ and a point $x_0 \in D$ the identity
\begin{equation}
u (x_0) \, a_m (\lambda r) = M (u, D) \label{new}
\end{equation}
holds for every function $u$ satisfying equation \eqref{Hh} in $D_r$. If also
\begin{equation}
D \subset B_{r_0} (x_0) , \ \ where \ \ \lambda r_0 = j_{m/2,1} \, ,
\label{exam}
\end{equation}
then $D = B_r (x_0)$.
\end{theorem}

\begin{remark}
{\rm For a given $\lambda$, this theorem is applicable only to domains, whose volume
is less than or equal to $|B_{r_0}|$, where $\lambda r_0 = j_{m/2,1}$, because every
such domain must lie within a ball of radius $r_0$.}
\end{remark}

Prior to proving this theorem, we introduce the following function:
\begin{equation} 
U (x) = a_{m-2} (\lambda |x|) \, , \quad x \in \RR^m . \label{U}
\end{equation}
Its main properties immediately follow from the representation:
\begin{equation}
U (x) = \frac{2 \, \Gamma (m/2)}{\sqrt \pi \, \Gamma ((m-1)/2)} \int_0^1 (1 -
s^2)^{(m-3)/2} \cos (\lambda |x| s) \, \D s \, , \label{PU}
\end{equation}
which is a consequence of Poisson's integral for $J_\nu$ (see \cite{NU}, p.~206).
The latter expression takes particularly simple form for $m=3$, namely, $U (x) =
(\lambda |x|)^{-1} \sin \lambda |x|$.

It is clear that this function solves \eqref{Hh} in $\RR^m$; indeed, it is easy to
differentiate \eqref{PU}, thus verifying the equation. Also, we have that $U (0) =
1$ and $U (x)$ decreases monotonically when $\lambda |x|$ belongs to the interval
$(0, j_{m/2,1})$, and is positive on the smaller interval $(0, j_{(m-2)/2,1})$.

\begin{proof}[Proof of Theorem 2.]
Without loss of generality, we suppose that the domain $D$ is located so that $x_0$
coincides with the origin. If we assume that $D \neq B_r (0)$, then $G_i = D
\setminus \overline{B_r (0)}$ and $G_e = B_r (0) \setminus \overline D$ are bounded
open sets such that $|G_e| = |G_i| \neq 0$ in view of the assumptions about $D$ and
$r$. In order to obtain a contradiction from the assumption made we write
identity~\eqref{new} for $U$ as follows:
\begin{equation}
|D| \, a_m (\lambda r) = \int_D U (y) \, \D y \, ; \label{1}
\end{equation}
here the condition $U (0) = 1$ is taken into account. Since property \eqref{MM}
holds for $U$ over $B_r (0)$, we write it in the same way:
\begin{equation}
|B_r| \, a_m (\lambda r) = \int_{B_r (0)} U (y) \, \D y \, . \label{2}
\end{equation}
Subtracting \eqref{2} from \eqref{1} and using the definition of $r$, we obtain
\begin{equation*}
0 = \int_{G_i} U (y) \, \D y - \int_{G_e} U (y) \, \D y < 0 \, .
\end{equation*}
Indeed, $U (y)$ monotonically decreases with $|y|$ in the whole $D$ because $D
\subset B_{r_0}$. Therefore, the difference is negative in view that $U (y)$ is
strictly greater than $[U (y)]_{|y| = r}$ in $G_e$ and strictly less than this value
in $G_i$, whereas $|G_i| = |G_e|$. The obtained contradiction proves the theorem.
\end{proof} 

Here, the argument is the same as in the proof of Theorem~1 given in \cite{Ku}; both
rely on monotonicity of a certain solution to the corresponding equation. However,
there is an essential distinction between the two theorems, and it concerns the size
of a domain. Indeed, no restriction on the size is imposed in Theorem~1, because its
proof is based on the function analogous to $U$, but increasing monotonically for
all $|x| > 0$. On the contrary, the radially symmetric function $U$ defined in
\eqref{U} decreases monotonically near the origin, but only when $\lambda |x|$
belongs to a bounded interval adjacent to zero. For this reason condition
\eqref{exam} is imposed in Theorem~2.

In the limit $\lambda \to 0$, equation \eqref{Hh} becomes the Laplace equation,
whereas $a_m (\lambda r) \to 1$, for which reason the assumption about $r$ becomes
superfluous in the limiting form of Theorem~2. Thus, letting $\lambda \to 0$ this
theorem turns into Kuran's.

In the formulation of Theorem~2, the integral $\int_D u (y) \, \D y$---it appears in
\eqref{new}---can be replaced by the flux integral $\int_{\partial D} \partial u /
\partial n_y \, \D S_y$ provided $\partial D$ is sufficiently smooth; here $n$ is
the exterior unit normal. Indeed, we have
\[ \int_D u (y) \, \D y = - \lambda^{-2} \int_D \nabla^2 u \, (y) \, \D y = - 
\lambda^{-2} \int_{\partial D} \partial u / \partial n_y \, \D S_y \, .
\]

To evaluate how restrictive is condition \eqref{exam}, we consider the following
simple example. Let a square membrane $D_s = \{ x \in \RR^2 : x_1, x_2 \in (0, a)
\}$, $a>0$, be fixed along the boundary. Then its free oscillations are described
by the following set of eigenfunctions:
\[ u_{ij} = \sin (i \pi x_1 / a) \, \sin (j \pi x_2 / a) \, , \quad u_{ji} = \sin (j 
\pi x_1 / a) \, \sin (i \pi x_2 / a) \, , \quad i,j = 1,2,\dots .
\]
They are linearly independent when $i \neq j$ and satisfy equation \eqref{Hh} with
\[ \lambda_{ij} = \lambda_{ji} = (\pi / a) \sqrt{i^2 + j^2} \, , \quad i,j = 1,2,\dots .
\] 
Let us check condition \eqref{exam} for $D_s$ with $\lambda = \lambda_{21} =
\lambda_{12} = \pi \sqrt 5 / a$. The reason to consider this particular value is
that $u_{21} (a/2, a/2) = u_{12} (a/2, a/2) = 0$ and
\[ \int_{D_s} u_{21} (x) \, \D x = \int_{D_s} u_{12} (x) \, \D x = 0 \, ,
\]
that is, identity \eqref{new} is valid for these functions provided $x_0 = (a/2,
a/2)$. It is clear that Theorem~2 is violated in this case, but what is the
difference between $\lambda r_0$ and $j_{1,1}$?

Since $r_0 = a / \sqrt 2$, we have that $\lambda_{21} r_0 = \pi \sqrt {5 / 2}
\approx 4.967294$ for $D_s$. Of course, the right-hand side in equality \eqref{exam}
is less than this value, namely, $j_{1,1} \approx 3.831706$. However, the difference
$\lambda r_0 - j_{1,1}$ is not too large.

In conclusion, we compare Theorem~2 with the result of Ramm \cite{R1,R3}, concerning
the so-called refined Schiffer's conjecture; the latter is similar to the celebrated
Serrin's theorem \cite{S}, but involves equation \eqref{Hh} instead of Poisson's.
Berenstein \cite{B}, p.~143, investigated this conjecture for simply connected
two-dimensional domains with smooth boundary (see also \cite{R2} for another
approach), whereas the original Schiffer's conjecture described in \cite{D} is
discussed in the review~\cite{CCH}.

In his monograph \cite{R1}, Ramm proved the following assertion, but it is more
convenient to give its formulation that appeared in the brief note \cite{R3}.

\begin{theorem}
Let $D \subset \RR^3$ be a bounded $C^2$-domain. If a nontrivial $u \in C^2 (D)
\cap C^1 (\overline D)$ satisfies equation \eqref{Hh} and the boundary conditions $u
= 0$, $\partial u / \partial n = c$ on $\partial D$ with constant~$c$, then $D$ is a
ball of radius $r$ such that $\lambda r$ is a zero of $j_0$--- the spherical Bessel
function of order zero.
\end{theorem}

In order to compare Theorems~3 and 2, we fix $\lambda > 0$ and notice that $j_0 (t)
= t^{-1} \sin t$. Therefore, Theorem~3 yields that $r = \pi k / \lambda$, $k =
1,2,\dots$, and so the existence of a solution to the overdetermined boundary value
problem, which is assumed in this theorem, guarantees that a domain is a ball
provided its volume belongs to a discrete sequence of values tending to infinity.

On the other hand, equality \eqref{new}---the key point of Theorem~2---implies that
a three-dimensional domain is a ball provided it lies within another ball of the
radius $r_0$ defined by the equality $\lambda r_0 = j_{3/2,1} \approx 4.493409$.
Therefore, the volume of a domain under consideration can be arbitrary, but within
the interval~$(0, |B_{r_0}|]$.

These considerations demonstrate clearly that the range of applicability of
Theorem~2 is completely different form that of Theorem~3.

\end{document}